\documentclass[12pt]{amsart}
\usepackage{amssymb}

\newtheorem{thm}{Theorem}[section]
\newtheorem{lem}[thm]{Lemma}
\newtheorem{cor}[thm]{Corollary}
\newtheorem{prop}[thm]{Proposition}
\newtheorem{ex}[thm]{Example}

\newtheorem*{prob*}{Open problem}

\theoremstyle{definition}

\newtheorem{defi}[thm]{Definition}

\theoremstyle{remark}

\newtheorem{rem}[thm]{Remark}
\newtheorem*{rem*}{Remark}

\DeclareMathOperator{\id}{id}

\DeclareMathOperator{\Hom}{Hom}

\newcommand{\kringel}{\mathbin{\raise1pt\hbox{$\scriptstyle\circ$}}}
\newcommand{\pkt}{\mathbin{\raise0pt\hbox{$\scriptstyle\bullet$}}}

\newcommand{\C}{\mathbb{C}}

\newcommand{\tr}{\mathop{\rm tr}}
\newcommand{\ad}{\mathop{\rm ad}}

\newcommand{\Ann}{\mathop{\rm Ann}}
\newcommand{\End}{\mathop{\rm End}}
\newcommand{\Der}{\mathop{\rm Der}}
\newcommand{\GDer}{\mathop{\rm GDer}}
\newcommand{\QDer}{\mathop{\rm QDer}}
\newcommand{\diag}{\mathop{\rm diag}}
\newcommand{\nil}{\mathop{\rm nil}}

\newcommand{\La}{\mathfrak{a}}
\newcommand{\Lb}{\mathfrak{b}}
\newcommand{\Lc}{\mathfrak{c}}

\newcommand{\Lg}{\mathfrak{g}}
\newcommand{\Lh}{\mathfrak{h}}

\newcommand{\Ll}{\mathfrak{l}}
\newcommand{\Ln}{\mathfrak{n}}

\newcommand{\Lr}{\mathfrak{r}}
\newcommand{\Ls}{\mathfrak{s}}

\newcommand{\CD}{\mathcal{D}}

\newcommand{\al}{\alpha}
\newcommand{\be}{\beta}
\newcommand{\ga}{\gamma}
\newcommand{\de}{\delta}
\newcommand{\ep}{\varepsilon}

\newcommand{\la}{\lambda}

\newcommand{\ov}{\overline}

\newcommand{\ra}{\rightarrow}

\renewcommand{\phi}{\varphi}

\begin{document}


\title[Post-Lie Algebra Structures]{Post-Lie algebra structures and generalized derivations
of semisimple Lie algebras}

\author[D. Burde]{Dietrich Burde}
\author[K. Dekimpe]{Karel Dekimpe}
\address{Fakult\"at f\"ur Mathematik\\
Universit\"at Wien\\
  Nordbergstr. 15\\
  1090 Wien \\
  Austria}
\email{dietrich.burde@univie.ac.at}
\address{Katholieke Universiteit Leuven\\
Campus Kortrijk\\
8500 Kortrijk\\
Belgium}
\email{karel.dekimpe@kuleuven-kortrijk.be}

\date{\today}

\thanks{The first author was supported by the FWF, Projekt P21683.}
\thanks{The second author expresses his gratitude towards the Erwin Schr\"odinger International
Institute for Mathematical Physics}

\begin{abstract}
We study post-Lie algebra structures on pairs of Lie algebras $(\Lg,\Ln)$, and 
prove existence results for the case that one of the Lie algebras 
is semisimple. For semisimple $\Lg$ and solvable $\Ln$ we show that there exist
no post-Lie algebra structures on  $(\Lg,\Ln)$. For semisimple $\Ln$ and certain solvable $\Lg$ we 
construct canonical post-Lie algebra structures. On the other hand we prove that there are no
post-Lie algebra structures for semisimple $\Ln$ and solvable, unimodular $\Lg$.
We also determine the generalized $(\al,\be,\ga)$-derivations of $\Ln$ in the semisimple case. 
As an application we classify post-Lie algebra structures induced by generalized derivations.
\end{abstract}

\maketitle

\section{Introduction}

Post-Lie algebras and post-Lie algebra structures recently have been introduced
in connection with homology of partition posets and the study of Koszul operads 
\cite{VAL}, \cite{LOD}. Surprisingly they also appeared in a quite different context
as well, namely in connection with nil-affine actions of Lie groups
\cite{BU41}. Post-Lie algebra structures generalize both LR-structures and 
pre-Lie algebra structures, which play a role in many areas, see
\cite{BU24},\cite{BU22},\cite{BU34},\cite{BU38}. In fact, a post-Lie algebra
structure on  $(\Lg,\Ln)$ with abelian Lie algebra $\Ln$ is just a pre-Lie algebra
structure on $\Lg$, and  a post-Lie algebra structure on $(\Lg,\Ln)$ with abelian 
Lie algebra $\Lg$ is an LR-structure on $\Ln$. In these extreme cases it is well known
that there are no such pre-Lie or LR-structures for semisimple Lie algebras.
Hence a natural question is whether there exist post-Lie algebra structures on
$(\Lg,\Ln)$ if one of the two Lie algebras is semisimple. Certainly there are
such structures if both $\Lg$ and $\Ln$ are semisimple, e.g., if they are isomorphic.
Therefore one might ask, whether it is possible that a pair $(\Lg,\Ln)$ admits
a post-Lie structure, if one of the Lie algebras is semisimple, and the other is
solvable. If $\Lg$ is semisimple, and $\Ln$ is solvable, it turns out that it is not
possible. A proof is given in section $4$. It uses the fact that a complex semisimple
Lie algebra $\Lg$ does not admit a pre-Lie algebra structure.
Conversely, if $\Ln$ is semisimple we construct canonical post-Lie algebra structures
on $(\Lg,\Ln)$ for some solvable, non-nilpotent Lie algebra $\Lg$. We consider the case
where $\Ln$ is semisimple and $\Lg$ is solvable and unimodular, i.e., all adjoint
operators have trace zero. We show that in this case there exist no post-Lie algebra 
structures. \\
In section $5$ we study generalized Lie algebra derivations. A particular type,
$(\al,\be,\ga)$-derivations, has been studied in connection with degeneration theory
of algebras, see \cite{BU36},\cite{NOH}. We determine the spaces of $(\al,\be,\ga)$-derivations 
for simple Lie algebras and study post-Lie algebra structures with semisimple $\Ln$, induced by
such generalized derivations.

\section{Post-Lie algebra structures}

Let $(\Lg, [\, ,])$ and $(\Ln, \{\, ,\})$ be two Lie brackets on the same vector space $V$ over
a field $k$. We call $(\Lg,\Ln)$ a {\it pair} of Lie algebras over $k$.
In particular we have $\dim (\Lg)=\dim (\Ln)$.
We define a post-Lie algebra structure as follows, see \cite{BU41}:

\begin{defi}\label{postlie}
Let $(\Lg, [\, ,])$ and $(\Ln, \{\, ,\})$ be two Lie brackets on a vector space $V$.
A {\it post-Lie algebra structure} on the pair $(\Lg,\Ln)$ is a $k$-bilinear product 
$x\cdot y$ satisfying the identities: 
\begin{align}
x\cdot y -y\cdot x & = [x,y]-\{x,y\} \label{post1}\\
[x,y]\cdot z & = x\cdot (y\cdot z) -y\cdot (x\cdot z) \label{post2}\\
x\cdot \{y,z\} & = \{x\cdot y,z\}+\{y,x\cdot z\} \label{post3}
\end{align}
for all $x,y,z \in V$. 
\end{defi}

The left multiplications of the algebra $A=(V,\cdot)$ are denoted by $L(x)$, i.e.,
we have $L(x)(y)=x\cdot y$ for all $x,y\in V$. Similarly the right multiplications
$R(x)$ are given by $R(x)(y)=y\cdot x$ for all $x,y\in V$. We note two immediate consequences of
the axioms \cite{BU41}:

\begin{lem}\label{2.2}
The map $L\colon \Lg\ra \End(V)$ given by $x\mapsto L(x)$ is a linear representation of the 
Lie algebra $\Lg$. Furthermore all operators $L(x)$ are Lie algebra derivations of $\Ln$.
\end{lem}

\begin{lem}
Let $x\cdot y$ be a post-Lie algebra structure on $(\Lg,\Ln)$. Then we have the following
identities:
\begin{align}
x\cdot \{y,z\}+y\cdot\{z,x\}+z\cdot\{x,y\} & = \{[x,y],z\}+ \{[y,z],x\}
+\{[z,x],y\} \label{post4}\\[0.2cm]
\{x,y\}\cdot z+\{y,z\}\cdot x+\{z,x\}\cdot y & = \{[x,y],z\}+ \{[y,z],x\}+ \{[z,x],y\} 
\label{post5}\\
 & \, +[\{x,y\},z]+[\{y,z\},x]+[\{z,x\},y]\nonumber
\end{align}
for all $x,y,z \in V$.
\end{lem}

\begin{ex}
If $\Ln$ is abelian, i.e., if $\{x,y\}=0$ for all $x,y\in V$,
then the conditions reduce to
\begin{align*}
x\cdot y-y\cdot x & = [x,y], \\
[x,y]\cdot z & = x\cdot (y\cdot z)-y\cdot (x\cdot z), 
\end{align*}
i.e., $x\cdot y$ is a {\it pre-Lie algebra structure} on the Lie algebra $\Lg$. 
\end{ex}

There is a large literature on pre-Lie algebras, also called left-symmetric algebras,
see \cite{BU24} for a survey. It is well known that a semisimple Lie algebra does not 
admit a pre-Lie algebra structure. It will be useful to repeat the proof.

\begin{prop}\label{2.5}
Let $\Lg$ be a semisimple Lie algebra over a field $k$ of characteristic zero.
Then $\Lg$ does not admit a pre-Lie algebra structure.
\end{prop}

\begin{proof}
We assume that $\Lg$ admits a pre-Lie algebra structure $x\cdot y$.
By lemma $\ref{2.2}$ we have $[L(x),L(y)]=L([x,y])$ for all $x,y\in \Lg$. 
Since $\Lg$ is semisimple we also have $[\Lg,\Lg]=\Lg$. These two conditions imply
that $\tr(L(x))=0$ for all $x\in \Lg$. In the same way we obtain $\tr(\ad(x))=0$ for 
the adjoint operators. The condition $x\cdot y-y\cdot x = [x,y]$ says that
$\ad (x)=L(x)-R(x)$ so that $\tr(R(x))=0$ for all $x\in \Lg$. It also says that
$\id \in Z^1(\Lg,\Lg_L)$, where $\Lg_L$ denotes the $\Lg$-module given by $L$.
By Whitehead's lemma, $\id$ is a $1$-coboundary. Hence there exists an element
$e\in \Lg$ with $R(e)=\id$. Taking traces we obtain that the identity map has trace zero.
Since the field $k$ has characteristic zero, and $\Lg\neq 0$, this is a contradiction.
\end{proof}

\begin{ex}
If $\Lg$ is abelian, then the conditions reduce to
\begin{align*}
x\cdot y-y\cdot x & = -\{x,y\} \\
x\cdot (y\cdot z)& = y\cdot (x\cdot z), \\
(x\cdot y)\cdot z & =(x\cdot z)\cdot y,
\end{align*}
i.e., $-x\cdot y$ is an {\it LR-structure} on the Lie algebra $\Ln$.
\end{ex}

For results on LR-structures on Lie algebras see \cite{BU33}, \cite{BU34}, \cite{BU38}. 
In particular we know that a semisimple Lie algebra does not admit an LR-structure, see 
\cite{BU34}.

\begin{prop}
Let $\Ln$ be a Lie algebra over a field of characteristic zero admitting an LR-structure.
Then $\Ln$ is solvable of class at most $2$.
\end{prop}

Recall that we have a correspondence between post-Lie algebra
structures on $(\Lg,\Ln)$ and embeddings $\Lg \hookrightarrow \Ln\rtimes \Der(\Ln)$, see \cite{BU41}.

\begin{prop}\label{2.8}
Let $x\cdot y$ be a post-Lie algebra structure on $(\Lg,\Ln)$. Then the map
\[
\phi\colon \Lg \ra \Ln\rtimes \Der(\Ln),\; x\mapsto (x,L(x))
\]
is an injective homomorphism of Lie algebras. Conversely any such embedding, with the identity
map on the first factor yields a post-Lie algebra structure on $(\Lg,\Ln)$.
\end{prop}

Here the bracket of $\Ln\rtimes \Der(\Ln)$ is given by  
\[
[(x,D),(x',D')]=(\{x,x'\}+D(x')-D'(x),[D,D']).
\]
Suppose that $x\cdot y$ is a post-Lie algebra structure on the pair of Lie algebras
$(\Lg,\Ln)$ such that $\Ln$ is centerless and satisfies $\Der(\Ln)=\ad (\Ln)$.
Then there is a $\phi\in \End(V)$ such that $x\cdot y=\{\phi(x),y \}$ for all $x,y\in V$, see
\cite{BU41}. This means $L(x)=\ad (\phi(x))$ for the linear operators $L(x)$.
In this case the axioms of a post-Lie algebra structure can be formulated as follows,
see \cite{BU41}:

\begin{prop}
Let $\Ln$ be a semisimple Lie algebra and  $\phi\in \End(V)$. Then
the product $x\cdot y=\{\phi(x),y \}$ is a post-Lie algebra structure on $(\Lg,\Ln)$
if and only if
\begin{align}
\{\phi(x),y\}+\{x,\phi(y)\} & =[x,y]-\{x,y\},\label{post6}\\
\phi([x,y]) & = \{\phi(x),\phi(y)\} \label{post7}
\end{align}
for all $x,y\in V$. 
\end{prop}
The second condition says that $\phi\colon \Lg\ra \Ln$ is a Lie algebra homomorphism. 

\begin{ex}
Let $\Ln$ be a semisimple Lie algebra. Then there are two obvious post-Lie algebra
structures on $(\Lg,\Ln)$ given by $\phi=0$ or $\phi=-\id$. In both cases
$\Lg$ is isomorphic to $\Ln$.
\end{ex}

Indeed, if $\phi=0$, the product $x\cdot y$ is zero, and $[x,y]=\{x,y\}$.
In the second case, where $\phi=-\id$, we have $x\cdot y=[x,y]=-\{x,y\}$. \\[0.2cm]
If $\Lg$ and $\Ln$ are both simple, these are the only possibilities.
We recall the following result \cite{BU41}:

\begin{prop}
Let $x\cdot y=\{\phi(x),y \}$ be a post-Lie algebra structure on $(\Lg,\Ln)$, where
$\Ln$ and $\Lg$ are simple. Then either $\phi=0$ or $\phi=-\id$. In both cases, 
$\Lg$ and $\Ln$ are isomorphic.
\end{prop}

The case where $\Lg$ and $\Ln$ are both semisimple is much more interesting.
Suppose that $\Lg$ and $\Ln$ are isomorphic. Then there exist the obvious post-Lie structures
arising from the simple factors. On the other hand, one can find much more post-Lie algebra
structures. We consider the following example. Let $\Ln=\Ls\Ll_2(\C)\oplus \Ls\Ll_2(\C)$, with
basis $e_1,f_1,h_1,e_2,f_2,h_2$ and Lie brackets
$$
\begin{array}{ll}
\{e_1,f_1\} = h_1, & \{e_2,f_2\} = h_2, \\
\{e_1,h_1\} = -2e_1, & \{e_2,h_2\} = -2e_2, \\
\{f_1,h_1\} = 2f_1, & \{f_2,h_2\} = 2f_2. 
\end{array}
$$
Now define $\phi=\begin{pmatrix} 0 & \vrule & 0 \\ \hline \\[-0.48cm] A & \vrule & 0
\end{pmatrix}$ by 
\[
A=\begin{pmatrix}
4  & -1 & -4 \\
-1 & 1 &  2 \\
-2 & 1 & 3
\end{pmatrix}.
\]

\begin{ex}\label{e210}
The product $x\cdot y=\{\phi(x),y \}$ with $\phi$ given as above 
defines a  post-Lie algebra structure on $(\Lg,\Ln)$,
where $\Lg$ and $\Ln$ are isomorphic to $\Ls\Ll_2(\C)\oplus \Ls\Ll_2(\C)$.
\end{ex}

The product is given by
$$
\begin{array}{lll}
e_1\cdot e_2 = -4e_2+h_2,  & f_1\cdot e_2 = 2e_2-h_2,  & h_1\cdot e_2 = 6e_2-2h_2, \\
e_1\cdot f_2 = 4f_2+4h_2,  & f_1\cdot f_2 = -2f_2-h_2, & h_1\cdot f_2 = -6f_2-4h_2, \\
e_1\cdot h_2 = -8e_2-2f_2, & f_1\cdot h_2 = 2e_2+2f_2, & h_1\cdot h_2 = 8e_2+4f_2.
\end{array}
$$
An easy calculation shows that this defines a post-Lie algebra structure on
$(\Lg,\Ln)$, where $\Lg$ is given by
$$
\begin{array}{lll}
[e_1,f_1]=h_1,                   & [f_1,h_1]=2f_1,      & [h_1,f_2]=-6f_2-4h_2, \\
\left[e_1,h_1\right]=-2e_1,      & [f_1,e_2]=2e_2-h_2,  & [h_1,h_2]=8e_2+4f_2, \\
\left[e_1,e_2\right]=-4e_2+h_2,  & [f_1,f_2]=-2f_2-h_2, & [e_2,f_2]=h_2, \\
\left[e_1,f_2\right]=4f_2+4h_2,  & [f_1,h_2]=2e_2+2f_2, & [e_2,h_2]=-2e_2, \\
\left[e_1,h_2\right]=-8e_2-2f_2, & [h_1,e_2]=6e_2-2h_2, &  [f_2,h_2]=2f_2.
\end{array}
$$
Furthermore, $\Lg$ is isomorphic to $\Ls\Ll_2(\C)\oplus \Ls\Ll_2(\C)$.

\begin{rem}
It is possible to compute all post-Lie algebra structures on
$(\Lg,\Ln)$, where $\Ln$ is $\Ls\Ll_2(\C)\oplus \Ls\Ll_2(\C)$ and
$\Lg$ is unimodular. The result is a list of products depending on parameters which
have to satisfy certain polynomial equations. It turns out that $\Lg$ is isomorphic 
to $\Ln$ in all cases. The matrix $\phi$ has one of the following forms:
\[
\begin{pmatrix} 0 & \vrule & 0 \\ \hline \\[-0.48cm] A & \vrule & 0 \end{pmatrix},\;
\begin{pmatrix} -\id & \vrule & 0 \\ \hline \\[-0.48cm] A & \vrule & 0 \end{pmatrix},\;
\begin{pmatrix} 0 & \vrule & 0 \\ \hline \\[-0.48cm] A & \vrule & -\id \end{pmatrix},\;
\begin{pmatrix} -\id & \vrule & 0 \\ \hline \\[-0.48cm] A & \vrule & -\id \end{pmatrix},\;
\]
and the transposed block types. \\
We give an example which generalizes $\ref{e210}$ as follows.
Let $\phi=\begin{pmatrix} 0 & \vrule & 0 \\ \hline \\[-0.48cm] A & \vrule & 0
\end{pmatrix}$ with 
\[
A=\begin{pmatrix}
\al & -\frac{\be^2}{4\al} & \be \\[0.2cm]
\ga &  -\frac{\de^2}{4\ga}  &  \de \\[0.2cm]
-\frac{\ep}{2} & -\frac{\be\de}{2\ep} & 1-\frac{\be\de}{2\al}
\end{pmatrix},
\]
where the parameters $\al,\be,\ga,\de,\ep$ have to satisfy $\al\ga\neq 0$, 
$\al\de-\be\ga\neq 0$, and the polynomial conditions $\ep=\al\de-\be\ga$, 
$\ep^2+4\al\ga=0$. For the choice
$(\al,\be,\ga,\de,\ep)=(4,-4,-1,2,4)$ we obtain example $\ref{e210}$.
\end{rem}

Finally we give an easy way to construct post-Lie algebra structures on $(\Lg,\Ln)$.
We use a direct vector space decomposition of $\Ln$ by subalgebras $\La$ and $\Lb$
as follows. 

\begin{prop}\label{2.9}
Let $(\Ln, \{\, ,\})$ be a Lie algebra which is a direct vector space sum  
$\Ln=\La\oplus \Lb$ of two subalgebras  $\La$ and $\Lb$. Let  $(\Lg, [\, ,])$ be the
Lie algebra given by $\La\oplus \Lb$ with the bracket
\begin{align*}
[a+b,a'+b'] & = \{a,a'\}-\{b,b'\}.
\end{align*}
Then we obtain a  post-Lie algebra structure on $(\Lg,\Ln)$ by 
\begin{align*}
(a+b)\cdot (a'+b') & = -\{b,a'+b'\}.
\end{align*}
\end{prop}

\begin{proof}
If we view the vector space $\La\oplus \Lb$ as a Lie algebra given by the direct Lie algebra sum 
of $\La$ and $\Lb$, then the map $\psi\colon \Lg\ra \La\oplus \Lb$, $\psi(a+b)=a-b$ is an 
algebra isomorphism. This shows that  $(\Lg, [\, ,])$ is indeed a Lie algebra. It is
different from $\Ln$ in general, because  $\Ln$ is not necessarily a direct Lie algebra sum 
of $\La$ and $\Lb$. \\
Let $x=a+b$, $y=a'+b'$ in $\La\oplus \Lb$. We have
\begin{align*}
x\cdot y -y\cdot x & = (a+b)\cdot (a'+b')-(a'+b')\cdot (a+b) \\
 & =  -\{b,a'+b'\} +\{b',a+b\} \\
 & = \{a',b\}+2\{b',b\}+\{b',a\} \\
 & = \{a,a'\}-\{b,b'\} - \{a+b,a'+b'\} \\
 & = [a+b,a'+b']-\{a+b,a'+b'\} \\
 & = [x,y]-\{x,y\}.
\end{align*}
This shows axiom \eqref{post1}. Furthermore we have
\begin{align*}
[x,y]\cdot z & = (\{a,a'\}-\{b,b'\})\cdot z\\
 & =  \{\{b,b'\},z\} \\
 & = \{b,\{b',z\}\}-\{b',\{b,z\}\} \\
 & = -\{b,b'\cdot z\}+\{b',b\cdot z\} \\
 & = (a+b)\cdot ((a'+b')\cdot z)-(a'+b')\cdot ((a+b)\cdot z) \\
 & = x\cdot (y\cdot z)-y\cdot (x\cdot z).
\end{align*}
This shows axiom \eqref{post2}. Finally we have
\begin{align*}
x\cdot \{y,z\} & = -\{b,\{y,z\}\}\\
 & = -\{\{z,b\},y\}-\{\{b,y\},z\} \\
 & = \{y,-\{b,z\}\}+ \{-\{b,y\},z\}\\
 & = \{y,x\cdot z\}+\{x\cdot y,z\},
\end{align*}
which shows \eqref{post3}.
\end{proof}

\section{Pairs $(\Lg,\Ln)$ with semisimple $\Ln$}

In this section we study post-Lie algebra structures on pairs $(\Lg,\Ln)$, where
$\Ln$ is semisimple. First we will use Proposition $\ref{2.9}$ to show
that there exist canonical post-Lie algebra structures such pairs 
$(\Lg,\Ln)$ where $\Ln$ is semisimple, or even simple, and $\Lg$ is some
solvable Lie algebra. \\
Suppose that $\Ln$ is a complex, semisimple Lie algebra with Cartan subalgebra
$\Lh$ and a root system of $\Ln$ with respect to  $\Lh$. If $\Ln^{+}$ denotes the
sum of the positive root spaces of $\Lg$ and $\Ln^{-}$ denotes the sum of the negative root spaces 
of $\Lg$, we obtain the so called {\it triangular decomposition} of $\Ln$,
\[
\Ln=\Ln^{-}\oplus \Lh\oplus \Ln^{+},
\]
which is a vector space direct sum of Lie subalgebras, where $\Lb^{+}=\Lh\oplus \Ln^{+}$ and
$\Lb^{-}=\Lh\oplus \Ln^{-}$ are Borel subalgebras of $\Ln$, and $\Ln^{\pm}$ is nilpotent. \\
As an example, a triangular decomposition of $\Ln=\Ls\Ll_n(\C)$ is given by
\begin{align*}
\Ls\Ll_n(\C) & = \Ln^{-}\oplus \Lh\oplus \Ln^{+} \\
 & = \bigoplus_{i>j}\C E_{ij}\oplus \bigoplus_{1\le i\le n-1} \C (E_{ii}-E_{i+1,i+1})\oplus 
\bigoplus_{i<j}\C E_{ij},
\end{align*}
where $E_{ij}$ is the $(n\times n)$-matrix with entry $1$ at position $(i,j)$ and zero entries 
otherwise.
Here $\Ln^{-}$ resp. $\Ln^{+}$ consists of strictly lower-triangular (resp. upper-triangular) matrices.

\begin{prop}\label{2.10}
Let $\Ln$ be a semisimple Lie algebra. Then there is a post-Lie algebra structure on
$(\Lg,\Ln)$ for some solvable non-nilpotent Lie algebra $\Lg$.
\end{prop}

\begin{proof}
Consider a triangular decomposition $\Ln=\Ln^{-}\oplus \Lh\oplus \Ln^{+}$ .
We can write the direct vector space decomposition $\Ln=\Ln^{-}\oplus \Lh\oplus \Ln^{+}$ as
$\Ln=\La\oplus\Lb$ such that one algebra is nilpotent and the other one is solvable and non-nilpotent.
We may choose $\La=\Ln^{-}$, $\Lb=\Lb^{+}$, or $\La=\Lb^{-}$, $\Lb=\Ln^{+}$, or interchange the
roles of $\La$ and $\Lb$. Then we obtain a post-Lie algebra structure on $(\Lg,\Ln)$ by 
Proposition $\ref{2.9}$. Here $\Lg$ is isomorphic to a
direct sum of a solvable and a nilpotent Lie algebra. Hence it is solvable. Since
$\Lb^{\pm}$ is non-nilpotent, $\Lg$ is non-nilpotent.
\end{proof}

We know from section $2$ that there exist many post-Lie algebra structures on pairs 
$(\Lg,\Ln)$, where both $\Ln$ and $\Lg$ are semisimple. It is perhaps less obvious that
there exist post-Lie algebra structures on pairs $(\Lg,\Ln)$, where $\Ln$ is semisimple, 
but $\Lg$ is not. The above Proposition shows that this is possible. If $\Ln$ is
semisimple, $\Lg$ can be solvable. 
On the other hand, $\Lg$ cannot be nilpotent in that
case because of the following result from \cite{BU41}.

\begin{prop}\label{n3.2}
Suppose that there is a post-Lie algebra structure on  $(\Lg,\Ln)$, where $\Lg$ 
is nilpotent. Then $\Ln$ must be solvable.
\end{prop}

Note that the solvable Lie algebra $\Lg$ we obtained in Proposition $\ref{2.10}$ is certainly not
nilpotent, since it is not even unimodular. Indeed, for semisimple $\Ln$ we can generalize
Proposition $\ref{n3.2}$ to the unimodular case.

\begin{thm}
Let $\Ln$ be a semisimple Lie algebra, and $\Lg$ be a solvable and unimodular
Lie algebra. Then there is no post-Lie algebra structure on $(\Lg,\Ln)$.
\end{thm}

\begin{proof}
Suppose there is a post-Lie algebra structure on  $(\Lg,\Ln)$. 
By \cite[Proposition 2.16]{BU41}, this means that we can view $\Lg$ as a subalgebra 
of $\Ln \times \Ln$, for which  the map 
$(p_1-p_2)_{\mid \Lg}\colon \Lg\ra \Ln$ is bijective. 
Here $p_1$ denotes the projection onto the first factor of $\Ln\times \Ln$, and 
$p_2$ the projection onto the second one. We remark here that in 
\cite{BU41}, we used the notation $\Ln \oplus  \Ln$ in stead of $\Ln \times \Ln$. However, in 
the sequel of this proof, we will also be dealing with internal direct sums and to avoid confusion, 
we will reserve the symbol $\oplus$ only for these internal sums. By assumption $p_1(\Lg)$ and 
$p_2(\Lg)$ are solvable, so that there exist Borel subalgebras $\Lb_1$ and $\Lb_2$ of $\Ln$ 
with $\Lb_1\supseteq p_1(\Lg)$ and $\Lb_2\supseteq p_2(\Lg)$. 
Now the intersection of two Borel subalgebras always contains a Cartan subalgebra $\Lh$,
and the dimension of all Borel subalgebras in $\Ln$ is equal to $(\dim(\Ln)+\dim (\Lh))/2$.
For a reference see Proposition $29.4.9$ in \cite{TAY}.
Hence let $\Lh \subseteq \Lb_1\cap \Lb_2 $ be a Cartan subalgebra of $\Ln$ of dimension $\ell$.
We have
\[
\dim (\Lb_1)=\dim (\Lb_2)=\frac{1}{2}(\dim (\Ln)+\ell).
\]
However, since $\Lb_1+\Lb_2 \supseteq p_1(\Lg)+p_2(\Lg)=\Ln$ and $\Lb_1,\Lb_2\subseteq \Ln$
we also have
\[
\dim (\Lb_1)+\dim (\Lb_2)-\dim (\Lb_1\cap \Lb_2)=\dim (\Ln).
\]
The two conditions imply $\dim (\Lb_1\cap \Lb_2)=\ell$. This means that
$\Lb_1\cap \Lb_2$ itself is a Cartan subalgebra of $\Ln$. 
Using the root space decomposition we obtain that $\Lb_1=\La_1\oplus \Lh$ and  
$\Lb_2=\La_2\oplus \Lh$ with nilpotent ideals $\La_1$ and $\La_2$. It follows that 
$\Ln=\La_1 \oplus \Lh \oplus \La_2$ and 
\[
\Lg\le (\La_1\oplus \Lh)\times  (\La_2\oplus \Lh) \le \Ln \times \Ln.
\]  
We will show that $\dim (\Lg\cap (\Lh \times  \Lh))=\ell$. To see this, write
$x=(a_1+h_1,a_2+h_2)\in \Lg$ with $a_i\in \La_i, h_i\in \Lh$ for $i=1,2$. Then
$(p_1-p_2)(x)=a_1-a_2+h_1-h_2 \in \Lh$ if and only if $x\in \Lh\times \Lh$.
In this case $a_1=a_2=0$ and $x=(h_1,h_2)\in \Lg\cap (\Lh \times \Lh)$.
It follows that the map
\[
(p_1-p_2)_{\mid \Lg\cap (\Lh \times \Lh)}\colon \Lg\cap (\Lh \times \Lh) \ra \Lh
\]
is bijective. Denote by $\ad_{\Ln}(x)$ the adjoint operators of $\Ln$, and by 
$\ad_{\Lg}(x)$ the adjoint operators of $\Lg$. For any $h\in \Lh$, we have that 
$\ad_{\Ln}(h)(\La_i)\subseteq  \La_i$ ($i=1,2$) and $\ad_{\Ln}(h) (\Lh)=0$. As $\Ln$ is semisimple,
$\tr \ad_{\Ln}(h)=0$ and hence we find that 
\begin{equation}\label{trace} 
0 = \tr {\ad}_{\Ln}(h)= \tr {\ad}_{\Ln}(h)_{\mid \La_1} + \tr {\ad}_{\Ln}(h)_{\mid \La_2}
\end{equation}
Now, consider any $x=(h_1,h_2)\in \Lg\cap  (\Lh \times \Lh)$ and let 
$y=(a_1+\tilde{h}_1, a_2 +\tilde{h}_2)$ be any element of $\Lg$. Then 
\[ {\ad}_{\Lg}(x)(y)= ({\ad}_{\Ln}(h_1) (a_1) , {\ad}_{\Ln}(h_2)(a_2) ).\]
From this (e.g.\ using explicit bases for the Lie algebras $\La_1$, $\La_2$ and $\Lh$) 
and the assumption that $\Lg$ is unimodular, it follows that 
\begin{align*}
0 & =\tr {\ad}_{\Lg}(x) =\tr {\ad}_{\Ln} (h_1)_{\mid \La_1}+\tr {\ad}_{\Ln}(h_2)_{\mid \La_2}\\
  & =\tr {\ad}_{\Ln}(h_1)_{\mid \La_1}-\tr {\ad}_{\Ln}(h_2)_{\mid \La_1}  \mbox{ \ \ (use (\ref{trace}))}\\
  & = \tr {\ad}_{\Ln} (h_1-h_2)_{\mid \La_1}.
\end{align*}
As the elements $h_1-h_2=(p_1-p_2)(x)$ are in 1--1 correspondence with the elements 
$h\in \Lh$, it follows that 
$\tr \ad_{\Ln} (h)_{\mid \La_1}=0$ for all $h\in \Lh$. We have $\Lb_1=\Lh\oplus \La_1$, so that  
$\tr \ad_{\Ln} (h)_{\mid \Lb_1}=0$ for all $h\in \Lh$. But then the Borel subalgebra $\Lb_1$
is unimodular. However, as the Borel subalgebra $\Lb_1$ can always be interpreted as being of the form 
$\Lb_1=\Lh\oplus \Ln^{+}$, where $\Ln^{+}$ is the sum of positive root spaces, there exist an $h\in \Lb_1$ 
for which $\ad_{\Ln} (h)_{\mid \Lb_1}$ only has non-negative (and at least one non-zero) 
real eigenvalues, which implies $\tr \ad_{\Ln} (h)_{\mid \Lb_1}>0$, a contradiction. 
\end{proof}
 
We close this section with some lemmas on eigenspaces for semisimple derivations of $\Ln$, which
we will need later. 

\begin{lem}\label{2.11}
Let $\Ln$ be a semisimple Lie algebra and $z\in \Ln$ an element such that
the adjoint operator $\ad(z)$ in $\Ln$ has three different eigenvalues
$0,\al$ and $\be$. Then $\al+\be=0$. 
\end{lem}

\begin{proof}
By looking at the semisimple part we may assume that  $\ad(z)$ is diagonalizable.
For a semisimple derivation $D\in \Der(\Ln)$ we have the eigenspace decomposition
$\Ln=\oplus V_{\la_i}$ corresponding to the different eigenvalues $\la_1,\ldots ,\la_n$ of $D$.
It holds $\{V_{\la_i},V_{\la_j}\}\subseteq V_{\la_i+\la_j}$. In particular, $\{V_{\la_i},V_{\la_j}\}=0$
if $\la_i+\la_j$ is not an eigenvalue of $D$. Now choose $D=\ad (z)$ and
denote by $V_0,V_{\al},V_{\be}$ the eigenspaces for the derivation $\ad(z)$ in $\Ln$
corresponding the the eigenvalues $0,\al$ and $\be$. By assumption we have
$\Ln=V_0\oplus V_{\al}\oplus V_{\be}$. Let $r,s\ge 1$ denote the multiplicities
of the eigenvalues $\al$ and $\be$. We have $\tr (\ad(z))=0$ since $\Ln$ is semisimple.
This implies $r\al+s\be=0$. Hence $\al+\al$ cannot be equal to $0$ or $\be$, because 
$r,s\ge 1$. Therefore $\al+\al$ is not an eigenvalue of $\ad(z)$, so that
$\{V_{\al},V_{\al}\}=0$. In the same way we obtain $\{V_{\be},V_{\be}\}=0$. Assume that $\al+\be \neq 0$.
Then $\{V_{\al},V_{\be}\}=0$ and $\{V_{\al},V_0\oplus V_{\al}\oplus V_{\be}\}\subseteq V_{\al}$. 
Hence $V_{\al}$ is an ideal of $\Ln$. We have $\ad(z)_{\mid V_{\al}}=\al\cdot \id$. On the other hand,
the trace of $\ad(z)$ restricted to an ideal of $\Ln$ is zero. Since $\al\neq 0$ this is a
contradiction. It follows $\al+\be=0$.
\end{proof}

Denote by $\al(\Ln)$ the maximal dimension of an abelian subalgebra of $\Ln$.

\begin{lem}\label{2.12}
Let $\Ln$ be a semisimple Lie algebra satisfying $\al(\Ln)\ge \frac{1}{3}\dim (\Ln)$.
Then equality holds and $\Ln$ is isomorphic to $\Ls\Ll_2(\C)\oplus \cdots \oplus  \Ls\Ll_2(\C)$,
where the number of summands is $\al(\Ln)= \frac{1}{3}\dim (\Ln)$.
\end{lem}

\begin{proof}
Let $\Ln=\Ls_1\oplus \cdots \oplus \Ls_k$ be the decomposition of $\Ln$ into simple
summands. Since the $\al$-invariant is additive it follows that
$\al(\Ln)=\al(\Ls_1)+\ldots +\al(\Ls_k)$.
For simple Lie algebras $\Ls_i$, the $\al$-invariant is well known, see \cite{SUT}. 
We always have $\al(\Ls_1)\le \frac{1}{3} \dim (\Ls_i)$, and equality
holds if and only if $\Ls_i$ is $\Ls\Ll_2(\C)$. Hence $\al(\Ln)\ge \frac{1}{3}\dim (\Ln)$
is only possible if $\al(\Ln)= \frac{1}{3}\dim (\Ln)$ and each summand $\Ls_i$ equals $\Ls\Ll_2(\C)$.
\end{proof}

\begin{lem}\label{2.13}
Let $\Ln$ be a semisimple Lie algebra which is a direct vector space sum of three
abelian subalgebras, i.e., $\Ln=\La\oplus \Lb\oplus \Lc$. Then  $\Ln$ is isomorphic to 
$\Ls\Ll_2(\C)\oplus \cdots \oplus  \Ls\Ll_2(\C)$.
\end{lem}

\begin{proof}
At least one of the three abelian subalgebras of $\Ln$ has a dimension bigger or equal than
$\frac{1}{3}\dim (\Ln)$. Therefore the $\al$-invariant of $\Ln$ satisfies
$\al(\Ln)\ge \frac{1}{3}\dim (\Ln)$. The claim follows from Lemma $\ref{2.12}$.
\end{proof}

\section{Pairs $(\Lg,\Ln)$ with semisimple $\Lg$}

We study here post-Lie algebra structures on pairs $(\Lg,\Ln)$, where
$\Lg$ is semisimple. It turns out that there exist no such structures
with a solvable Lie algebra $\Ln$. For pairs  $(\Lg,\Ln)$ with $\Lg$ semisimple and
$\Ln$ abelian this is Proposition $\ref{2.5}$. \\
Let $M$ be a $\Lg$-module and $m\in M$. Denote by $\Ann (m)=\{ x \in \Lg \mid x.m=0\}$
the annihilator of $m$ in $\Lg$.

\begin{lem}\label{3.1}
Let $\Lg$ be a semisimple Lie algebra and $M$ be a $\Lg$-module with $\dim (M)\le \dim (\Lg)$.
Then $\Ann (m)\neq 0$ for all $m\in M$.
\end{lem}

\begin{proof}
Consider the linear map $\psi\colon \Lg\ra M$ given by $x\mapsto x.m$. Its kernel is
$\Ann (m)$. If $\dim (M)< \dim (\Lg)$ then $\psi$ is not injective and has a non-trivial 
kernel. If $\dim (M)= \dim (\Lg)$, then $\psi$ is a $1$-cocycle in $Z^1(\Lg,M)$.
Assume that it is bijective. Then consider the $\Lg$-module $\ov{M}$ given by 
$\psi^{-1}\rho \psi$, where $\rho$ is the representation associated to $M$.
It is easy to see that we would obtain $\id \in Z^1(\Lg,\ov{M})$, defining a pre-Lie algebra 
structure on $\Lg$. This contradicts Proposition $\ref{2.5}$.
It follows that $\ker (\psi)=\Ann(m)$ is non-trivial.
\end{proof}

The lemma is no longer true if $\dim (M)>\dim (\Lg)$.

\begin{thm}\label{3.2}
Let $(\Lg,\Ln)$ be a pair of Lie algebras, where $\Lg$ is semisimple and $\Ln$ is solvable.
Then there is no post-Lie algebra structure on $(\Lg,\Ln)$.
\end{thm}

\begin{proof}
Assume that there is a post-Lie algebra structure on $(\Lg,\Ln)$. Then 
there exists an embedding $\Lg  \hookrightarrow \Ln\rtimes \Der(\Ln)$ for which
${p_1}_{\mid \Lg}\colon \Lg\ra \Ln$ is a bijection, see Proposition $\ref{2.8}$.
Here we identify $\Lg$ with its image in $\Ln\rtimes \Der(\Ln)$ and denote by
$p_1\colon  \Ln\rtimes \Der(\Ln) \ra \Ln$ the projection onto the first
factor. Let $p_2\colon  \Ln\rtimes \Der(\Ln) \ra \Der(\Ln)$ be the projection onto
the second factor and $\Ls=p_2(\Lg)$. We have $\Lg\simeq \Ls$ and 
$\Lg\le \Ln\rtimes \Ls$. Both $\Lg$ and $\Ls$ are Levi subalgebras of $\Ln\rtimes \Ls$. 
By the Levi-Malcev Theorem these are conjugated by an inner automorphism of the form
$\exp(\ad (z))$ for some $z\in \nil(\Ln)$ in the nilradical of $\Ln$.
Fix an element $w\in \Ln$ such that $\Lg=\exp(\ad (w))(\Ls)$. For $x\in \Ls$ we have
\[
\ad(w)(x)=\{w,x\}=-x.w,
\]
where the dot denotes the action of $\Ls$ on $\Ln$ in the semidirect product $\Ln\rtimes \Ls$.
As a vector space we may view $\Ln$ as an $\Ls$-module. By Lemma $\ref{3.1}$ there exists
a nonzero element $x\in \Ls$ such that $0=x.w$. It follows that $\ad(w)(x)=0$ and hence
$\ad(w)^n(x)=0$ for all $n\ge 1$. This implies that
\[
\exp(\ad(w))(x)=x
\]
for this nonzero $x \in \Ls$. It follows $x\in \Lg$ and hence $p_1(x)=0$.
This contradicts the fact that ${p_1}_{\mid \Lg}$ is bijective.  
\end{proof}

\section{Generalized derivations of Lie algebras}

In the next section we will construct a special kind of post-Lie algebra structures using
certain generalized Lie algebra derivations. Indeed, there are  various generalizations of them in the literature.
We will be concerned here with $\de$-derivations of Lie algebras \cite{FIL1}, \cite{FIL2},
\cite{ZUS}, with generalized derivations in the sense of Leger and Luks \cite{LEL}, and with
$(\al,\be,\ga)$-derivations of Lie algebras \cite{NOH},\cite{BU36}.
The classes of such derivations are partly overlapping, see the definitions below. 
We will present some results for generalized derivations of simple and semisimple Lie algebras.
In particular, we will determine all spaces of $(\al,\be,\ga)$-derivations for 
complex simple Lie algebras. \\[0.2cm]
Let $\Lg$ be a finite-dimensional Lie algebra over a field $k$. 
Denote by $\End(\Lg)$ the space of linear maps $\phi\colon \Lg\ra \Lg$. They are not
necessarily Lie algebra homomorphisms. 

\begin{defi}\cite{NOH}
For arbitrary scalars $\al,\be,\ga \in k$ we define the vector space $\CD(\al,\be,\ga)$
by
\[
\CD(\al,\be,\ga)=\{ \phi \in \End(\Lg)\mid \al\phi([x,y])=\be [\phi(x),y]+\ga [x,\phi(y)]
\;\forall \; x,y\in \Lg \}.
\]
These linear maps are called {\it $(\al,\be,\ga)$-derivations} of $\Lg$.
\end{defi}

Note that $\CD(1,1,1)=\Der(\Lg)$ is just the space of derivations of $\Lg$ in the ordinary sense.
Furthermore, many special cases have been considered before.
We list the definitions for the special cases $(\al,\be,\ga)=(1,\de,\de), (1,1,0), (0,1,-1)$.

\begin{defi}[Filippov]
For any $\de \in k$ define the vector space of {$\de$-derivations} by
\[
\CD(1,\de,\de)=\{  \phi \in \End(\Lg)\mid \phi([x,y])=\de [\phi(x),y]+\de [x,\phi(y)] \; \forall
\; x,y\in \Lg\}.
\]
\end{defi}

\begin{defi}
The space of $(1,1,0)$-derivations is given by
\begin{align*}
\CD(1,1,0) & = \{ \phi \in \End(\Lg)\mid \phi([x,y])= [\phi(x),y] \;\forall \; x,y\in \Lg\}. 
\end{align*}
Note that $\phi\in \CD(1,1,0)$ also satisfies $\phi([x,y])=-\phi ([y,x])=-[\phi(y),x]=[x,\phi(y)]$.
Hence
\begin{align*}
\CD(1,1,0) & =\{ \phi \in \End(\Lg) \mid \phi\circ \ad(x)=\ad(x)\circ \phi \;\forall\; x\in \Lg\}\\
           & = C(\Lg),
\end{align*}
which is called the {\it centroid} of $\Lg$.
\end{defi}

The centroid of $\Lg$ is an associative subalgebra of $\End(\Lg)$, such that
\[
[C(\Lg),C(\Lg)]\subseteq \Hom(\Lg/[\Lg,\Lg],Z(\Lg)).
\]
In particular, if $\Lg$ is perfect or centerless, $C(\Lg)$ is commutative.

\begin{defi}
The space of $(0,1,-1)$-derivations is given by
\begin{align*}
\CD(0,1,-1) & = \{ \phi \in \End(\Lg)\mid [\phi(x),y] = [x,\phi(y)] \;\forall \; x,y\in \Lg\} \\
            & = QC(\Lg),
\end{align*}
which is called the {\it quasicentroid} of $\Lg$.
\end{defi}

The quasicentroid is also a special case of the generalized Lie algebra derivations 
defined by Leger and Luks. The definition is as follows.

\begin{defi}[Leger, Luks]
The space of {\it generalized derivations} of $\Lg$ is given by
\[
\GDer(\Lg)=\{ \phi\in \End(\Lg)\mid \;\exists \;\sigma,\tau \in \End(\Lg) \text{ with } 
\tau([x,y])=[\phi(x),y]+[x,\sigma(y)] \;\forall\; x,y \in \Lg\}.
\]
\end{defi}

An important special case is the space of quasiderivations.

\begin{defi}
The space of {\it quasiderivations} of $\Lg$ is given by
\[
\QDer(\Lg)=\{ \phi\in \End(\Lg)\mid \;\exists \;\tau \in \End(\Lg) \text{ with } 
\tau([x,y])=[\phi(x),y]+[x,\phi(y)] \;\forall\; x,y \in \Lg\}.
\]
\end{defi}

We have the inclusions
\[
\ad(\Lg)\subseteq \Der(\Lg)\subseteq \QDer(\Lg)\subseteq \GDer(\Lg)\subseteq \End(\Lg),
\]
where $\ad(\Lg)$ denotes the space of inner deriavtions of $\Lg$. 
Furthermore we have
\begin{align*}
\QDer(\Lg)+QC(\Lg) & = \GDer(\Lg),\\
\Der(\Lg)+ C(\Lg) & \subseteq  \QDer(\Lg),
\end{align*}
where the first equality requires $k$ to be a field of characteristic
different from $2$.
Leger and Luks described conditions on $\Lg$ which force $QC(\Lg)=C(\Lg)$, or equivalently 
$\GDer(\Lg)=\QDer(\Lg)$. Theorem $5.28$ of \cite{LEL} proves the following result.
We assume that $k$ is the field of complex numbers.

\begin{thm}[Leger, Luks]\label{3.7}
Let $\Lg$ be a complex centerless and perfect Lie algebra. Then $QC(\Lg)=C(\Lg)$ and
$\GDer(\Lg)=\QDer(\Lg)$.
\end{thm}

Recall that the assumptions mean $Z(\Lg)=0$ and $[\Lg,\Lg]=\Lg$.

\begin{cor}
Let $\Lg$ be a complex simple Lie algebra. Then 
\[
\CD(0,1,-1)=QC(\Lg)=C(\Lg)=\C\cdot\id.
\]
\end{cor}

\begin{proof}
Indeed, by Schur's lemma we have $C(\Lg)=\C\cdot \id$. The claim follows 
immediately from Theorem $\ref{3.7}$.
\end{proof}

If $\Lg$ is simple we have $\Der(\Lg)\oplus C(\Lg)=\ad (\Lg)\oplus \C\cdot \id \subseteq
\QDer(\Lg)$. The next result shows that we have equality for simple Lie algebras different from
$\Ls\Ll_2(\C)$, see Corollary $4.16$ of \cite{LEL}.

\begin{thm}[Leger, Luks]\label{3.9}
Suppose that $\Lg$ is a complex simple Lie algebra of rank at least two. Then 
$\QDer(\Lg)=\ad (\Lg)\oplus \C\cdot \id$.
\end{thm}

In fact, for $\Ls\Ll_2(\C)$ the result is different.

\begin{prop}\label{3.10}
For $\Lg=\Ls\Ll_2(\C)$ we have $\QDer(\Lg)=\End(\Lg)$.
\end{prop}

\begin{proof}
Let $(e,f,h)$ be the standard basis of $\Ls\Ll_2(\C)$ with $[e,f]=h$, $[h,e]=2e$ 
and $[h,f]=-2f$. For a given $\phi\in \End(\Lg)$ with matrix $\phi=(x_{ij})$ define 
$\tau\in \End(\Lg)$ by
\[
\tau=\begin{pmatrix}
x_{11}+x_{33} & -x_{12}       & -2x_{32} \\
-x_{21}       & x_{22}+x_{33} & -2x_{31} \\
-\frac{x_{23}}{2}  & -\frac{x_{13}}{2}  & x_{11}+x_{22}
\end{pmatrix}.
\]
Then it is easy to see that $\tau$ satisfies
\[
\tau([a,b])=[\phi(a),b]+[a,\phi(b)]
\]
for all $a,b\in \Lg$. Hence every $\phi\in \End(\Lg)$ is a quasiderivation of $\Lg$.
\end{proof}

Note that the space $\ad (\Lg)\oplus \C\cdot \id$ is $4$-dimensional in this case,
whereas $\QDer(\Lg)$ is $9$-dimensional. \\[0.2cm]
We want to determine now the spaces $\CD(\al,\be,\ga)$ for complex Lie algebras. 
The following elementary result was shown in \cite{NOH}.

\begin{prop}
Let $\Lg$ be a complex Lie algebra and $\al,\be,\ga \in \C$. Then $\CD(\al,\be,\ga)$
equals one of the following subspaces of $\End(\Lg)$.
\begin{itemize}
\item[(a)]  $\CD(0,0,0)=\End(\Lg)$.
\item[(b)]  $\CD(1,0,0)=\{\phi \in \End(\Lg)\mid \phi([\Lg,\Lg])=0 \}$.
\item[(c)]  $\CD(0,1,-1)=QC(\Lg)$.
\item[(d)]  $\CD(1,1,-1)=\CD(0,1,-1)\cap \CD(1,0,0)$.
\item[(e)]  $\CD(\de,1,1)$, $\de \in \C$.
\item[(f)]  $\CD(\de,1,0)=\CD(0,1,-1)\cap \CD(2\de,1,1)$, $\de \in \C$.
\end{itemize}
\end{prop}

Note that $\CD(1,0,0)$ is a vector space of dimension $\dim (\Lg/[\Lg,\Lg])\cdot \dim (\Lg)$,
and 
\[
\CD(0,1,0)=\{\phi \in \End(\Lg)\mid \phi(\Lg)\subseteq Z(\Lg)\}
\]
is a vector space of dimension $\dim (Z(\Lg))\cdot \dim (\Lg)$.
The list shows that it is enough to determine the cases $(c)$ and $(e)$, because then
all other cases are determined as well.

\begin{cor}
Let $\Lg$ be a complex, centerless and perfect Lie algebra. Then $\CD(0,1,-1)=QC(\Lg)=C(\Lg)$
is a commutative, associative subalgebra of $\End(\Lg)$, and
\[
\CD(1,0,0)=\CD(1,1,-1)=\CD(0,1,0)=0.
\]
\end{cor}

\begin{proof}
By Theorem $\ref{3.7}$ we have $QC(\Lg)=C(\Lg)$. This is an associative subalgebra
of $\End(\Lg)$ with $[C(\Lg),C(\Lg)]\subseteq \Hom(\Lg/[\Lg,\Lg],Z(\Lg))$. If 
$\Lg$ is centerless or perfect, then this space is zero.
\end{proof}

It remains to consider the spaces $\CD(\de,1,1)$ for $\de \in \C$. 
For $\de=0$ we obtain the following result from \cite{LEL}, Lemma $6.1$.

\begin{prop}
Let $\Lg$ be a complex semisimple Lie algebra. Then $\CD(0,1,1)=0$.
\end{prop}

If $\de\neq 0$, then $\CD(\de,1,1)=\CD(1,\frac{1}{\de},\frac{1}{\de})$, and we deal with
$\frac{1}{\de}$-derivations in the sense of Filippov. We obtain the following result 
from \cite{FIL2}, Theorem $2$, and \cite{FIL1}, Theorem $6$.

\begin{thm}
Let $\Lg$ be a complex simple Lie algebra and $\de\neq 0,1,-1,2$. Then  $\CD(\de,1,1)=0$.
Furthermore we have  $\CD(2,1,1)=\C\cdot \id$ and $\CD(1,1,1)=\ad (\Lg)$. 
\end{thm}

For the case $\de=-1$ we obtain the following result from  \cite{FIL2}, Theorem $4$.

\begin{thm}
Let $\Lg$ be a complex simple Lie algebra of rank at least two. Then
\[
\CD(-1,1,1)=0.
\]
\end{thm}

As in the case of Theorem $\ref{3.9}$, also here the result does not hold for
$\Ls\Ll_2(\C)$.

\begin{ex}
For $\Lg=\Ls\Ll_2(\C)$ with respect to the standard basis, the space
$\CD(-1,1,1)$ is the linear span of
\[
\begin{pmatrix} 1 & 0 & 0 \\ 0 & 1 & 0 \\ 0 & 0 & -2 \end{pmatrix},\; 
\begin{pmatrix} 0 & 0 & 0 \\ 1 & 0 & 0 \\ 0 & 0 & 0 \end{pmatrix},\;
\begin{pmatrix} 0 & 1 & 0 \\ 0 & 0 & 0 \\ 0 & 0 & 0 \end{pmatrix},\;
\begin{pmatrix} 0 & 0 & 0 \\ 0 & 0 & 1 \\ 0 & 0 & 0 \end{pmatrix},\;
\begin{pmatrix} 0 & 0 & 2 \\ 0 & 0 & 0 \\ 0 & 1 & 0 \end{pmatrix}.
\]
\end{ex}

This is a complement of the vector space $\ad(\Lg)\oplus \C\cdot \id$ in $\End(\Lg)$,
see also Proposition $\ref{3.10}$.

\begin{rem}
It holds $\CD(\de,1,1)\subseteq \QDer(\Lg)$ for all $\de\in\C$, and we can derive
the previous results also from Theorem $\ref{3.9}$. Indeed, let $\phi\in \CD(\de,1,1)$,
and $\Lg$ be a simple Lie algebra of rank at least $2$.
Then there is a $z\in \Lg$ and a $\la\in \C$ such that $\phi(x)=[z,x]+\la x$. By assumption
we have
\[
\de\cdot \phi([x,y]) = [\phi(x),y]+[x,\phi(y)]
\]
for all $x,y\in \Lg$. If we substitute $\phi$ as above, this condition can be written as
\[
(\de-1)\cdot [z,[x,y]]=(2-\de)\la \cdot [x,y]
\]
for all $x,y\in \Lg$. Now an easy case distinction shows that  $\CD(\de,1,1)=0$ for all
$\de\neq 1,2$, and $\CD(2,1,1)=\C\cdot \id$,  $\CD(1,1,1)=\ad (\Lg)\simeq \Lg$.
\end{rem}

We summarize the results on $\CD(\al,\be,\ga)$ for simple Lie algebras.

\begin{prop}
Let $\Lg$ be a simple complex Lie algebra of dimension $\dim (\Lg)\ge 4$, and $\al,\be,\ga \in \C$. 
Then $\CD(\al,\be,\ga)$ equals one of the following subspaces of $\End(\Lg)$.
\begin{itemize}
\item[(a)]  $\CD(0,0,0)=\End(\Lg)$.
\item[(b)]  $\CD(1,0,0)=0$.
\item[(c)]  $\CD(0,1,-1)=\C\cdot \id$.
\item[(d)]  $\CD(1,1,-1)=0$.
\item[(e)]  $\CD(\de,1,1)=0$ for all $\de \neq 1,2$; \enskip $\CD(1,1,1)=\ad (\Lg)$; \enskip
 $\CD(2,1,1)=\C\cdot \id$. 
\item[(f)]  $\CD(\de,1,0)=0$ for all $\de \neq 1$; \enskip $\CD(1,1,0)=\C\cdot \id$.
\end{itemize}
\end{prop}

\section{Post-Lie algebra structures induced by generalized derivations}

Let $\Ln$ be semisimple. Suppose that $x\cdot y$ is a post-Lie algebra structure on the 
pair of Lie algebras $(\Lg,\Ln)$. Then there is a $\phi\in \End(V)$ such that 
$x\cdot y=\{\phi(x),y \}$ for all $x,y\in V$.
Recall that a product $x\cdot y=\{\phi(x),y \}$ with a linear map $\phi\colon \Lg\ra \Ln$
is a post-Lie algebra structure on $(\Lg,\Ln)$
if and only $\phi$ satisfies the conditions \eqref{post6} and \eqref{post7}.
Condition \eqref{post7} says that $\phi$ is a Lie algebra homomorphism, whereas
condition \eqref{post6} looks similar to a derivation condition for $\phi$, but involves two Lie 
brackets. 
However, if we assume that the Lie bracket of $\Lg$ is given as a linear function of the Lie
bracket of $\Ln$, then $\phi$ is indeed a quasiderivation of $\Ln$.

\begin{prop}
Let $x\cdot y$ be a post-Lie algebra structure on $(\Lg,\Ln)$ with 
$\Ln$ semisimple, and $x\cdot y=\{\phi(x),y \}$ for some $\phi\in \End(V)$.
Assume that $[x,y]=\tau (\{ x,y \})$ for some $\tau \in \End(V)$. 
Then $\phi\in \QDer(\Ln)$.
\end{prop}

\begin{proof}
Since $\Ln$ is semisimple, it is centerless and satisfies $\Der(\Ln)=\ad (\Ln)$.
Condition \eqref{post6} says
\[
\{\phi(x),y\}+\{x,\phi(y)\}=(\tau-\id)(\{x,y\}),
\]
which means $\phi\in \QDer(\Ln)$. 
\end{proof}

The space $\QDer(\Ln)$ has been studied. For any Lie algebra $\Ln$ we have
\[
\Der(\Ln)+ C(\Ln)\subseteq \QDer(\Ln).
\]
Suppose that $\Ln$ is simple of rank at least $2$. Then  $\QDer(\Ln)=\Der(\Ln)\oplus C(\Ln)=\ad(\Ln)\oplus
\C\cdot \id$ by Theorem $\ref{3.9}$. This yields the following result.

\begin{prop}\label{4.2}
Suppose that $x\cdot y$ is a post-Lie algebra structure on $(\Lg,\Ln)$, where
$\Ln$ is simple of rank at least $2$. Assume that $[x,y]=\tau (\{ x,y \})$ for some 
$\tau \in \End(V)$. Then 
\begin{align}\label{8}
x\cdot y & =\{\{z,x\},y\}+\la \{x,y\} 
\end{align}
for some $z\in \Ln$ and some $\la\in \C$.
\end{prop}

\begin{proof}
We write $x\cdot y=\{\phi(x),y \}$ for some $\phi\in \End(V)$. By the proposition we have
$\phi\in \QDer(\Ln)=\ad(\Ln)\oplus \C\cdot \id$. Hence there is a $z\in \Ln$ and a $\la\in\C$ such
that $\phi(x)=\{z,x\}+\la x$. This shows the claim.
\end{proof}

The assumption that $[x,y]=\tau (\{ x,y \})$ is of course a restriction on $(\Lg,\Ln)$.
Proposition $\ref{2.10}$ yields obvious examples of post-Lie structures on $(\Lg,\Ln)$
where $\Ln$ is simple, and the Lie brackets of $\Lg$ and $\Ln$ do not satisfy
a relation as above. Then also the conclusion of Proposition $\ref{4.2}$ need not hold. \\
We present an explicit example for $\Ln=\Ls\Ll_3(\C)$. Consider the following basis:
\begin{align*}
e_1 & =E_{12},\; e_2 = E_{13},\; e_3=E_{21},\;  e_4=E_{23},\;e_5 =E_{31}, \\
e_6 & = E_{32},\; e_7=E_{11}-E_{22},\;  e_8=E_{22}-E_{33}. 
\end{align*}
Then the Lie brackets are defined by
\begin{align*}
\{e_1,e_3 \} & = e_7,\;\{e_1,e_4 \} = e_2,\;\{e_1,e_5 \} = -e_6,\; \{e_1,e_7 \} = -2e_1,\; 
\{e_1,e_8 \} = e_1, \\
\{e_2,e_3 \} & = -e_4,\;\{e_2,e_5 \} = e_7+e_8,\;\{e_2,e_6 \} = e_1,\; \{e_2,e_7 \} = -e_2,\;
\{e_2,e_8 \} = -e_2, \\
\{e_3,e_6 \} & = -e_5,\;\{e_3,e_7 \} = 2e_3,\;\{e_3,e_8 \} = -e_3, \\
\{e_4,e_5 \} & = e_3,\; \{e_4,e_6 \} = e_8,\; \{e_4,e_7 \} = e_4,\;\{e_4,e_8 \} = -2e_4, \\
\{e_5,e_7 \} & = e_5,\;\{e_5,e_8 \} = e_5, \\
\{e_6,e_7 \} & = -e_6,\;\{e_6,e_8 \} = 2e_6. 
\end{align*}
Let $\Ls\Ll_3(\C)= \Ln^{-}\oplus \Lh\oplus  \Ln^{+}$ be the triangular decomposition
given by $\Ln^{-}=\langle e_3,e_5,e_6\rangle$, $\Lh=\langle e_7,e_8\rangle $ and
$\Ln^{+}=\langle e_1,e_2,e_4\rangle$.

\begin{ex}
The triangular decomposition $\Ls\Ll_3(\C)=\La\oplus \Lb$ with $\La=\Lh\oplus \Ln^{+}$ and
$\Lb=\Ln^{-}$ induces a post-Lie algebra structure on $(\Lg,\Ls\Ll_3(\C))$ which is not of the
form $x\cdot y =\{\{z,x\},y\}+\la \{x,y\}$ for some $z\in \Ln$ and some $\la\in \C$.
\end{ex}

By Proposition $\ref{2.9}$ the brackets of $\Lg$ are given by
\begin{align*}
[e_1,e_4] & = e_2,\;[e_1,e_7]= -2e_1,\;[e_1,e_8]= e_1, \\
[e_2,e_7 ] & = -e_2,\;[e_2,e_8]= -e_2,\;[e_3,e_6]= e_5, \\
[e_4,e_7 ] & = e_4,\; [e_4,e_8]= -2e_4. 
\end{align*}

This Lie algebra is $3$-step solvable and non-nilpotent. It is of dimension $8$ and  has a 
$1$-dimensional center.
Assume that $[x,y]=\tau (\{ x,y \})$ for some map $\tau\colon \Lg\ra \Ln$.
Then
\begin{align*}
e_1 & = [e_1,e_8] = \tau(\{e_1,e_8\})=\tau(e_1),\\
0 & = [e_2,e_6] = \tau(\{e_2,e_6\})=\tau(e_1),
\end{align*}
which is a contradiction. The post-Lie structure is given by $x\cdot y=\{\phi(x),y\}$ where
\[
\phi=\diag (0,0,-1,0,-1,-1,0,0).
\]
More precisely, the non-zero products are given by
\begin{align*}
e_3\cdot e_1 & = e_7,\; e_3\cdot e_2= -e_4,\;e_3\cdot e_6= e_5,\;e_3\cdot e_7= -2e_3,\; e_3\cdot e_8= e_3,\\
e_5\cdot e_1 & = -e_6,\; e_5\cdot e_2= e_7+e_8,\;e_5\cdot e_4= e_3,\;e_5\cdot e_7= -e_5,\; e_5\cdot e_8= -e_5,\\
e_6\cdot e_2 & = e_1,\; e_6\cdot e_3= -e_5,\;e_6\cdot e_4= e_8,\;e_6\cdot e_7= e_6,\; e_6\cdot e_8= -2e_6.
\end{align*}
It can be easily checked that $\phi$ is not of the form $\ad (z)+\la \cdot \id$. \\[0.2cm]
The result of Proposition $\ref{4.2}$ is a motivation to study post-Lie algebra structures 
on $(\Lg,\Ln)$  of the form \eqref{8}, where $\Ln$ is semisimple:

\begin{thm}\label{4.4}
Let $\Ln$ be a semisimple Lie algebra which is not a direct sum of copies of $\Ls\Ll_2(\C)$.
Suppose that
\[
x\cdot y =\{\{z,x\},y\}+\la \{x,y\} 
\]
is a post-Lie algebra structure on $(\Lg,\Ln)$ for some $z\in \Ln$, $\la\in \C$. 
Then either $x\cdot y=0$ and $[x,y]=\{x,y\}$, or  $x\cdot y=-\{x,y\}$ and $[x,y]=-\{x,y\}$.
\end{thm}

\begin{proof}
We write  $x\cdot y=\{\phi(x),y \}$ with $\phi(x)=\{z,x\}+\la x$. This is a post-Lie algebra
structure on $(\Lg,\Ln)$ if and only if the identities \eqref{post6} and \eqref{post7} are
satisfied, which means
\begin{align}
[x,y] & = \{z,\{x,y\}\}+(2\la+1)\{x,y\}, \label{post9}\\
\{\{z,x\},\{z,y\}\} & = \{z,\{z,\{x,y\}\}\}+(2\la+1)\{z,\{x,y\}\}+(\la^2+\la)\{x,y\} \label{post10}
\end{align}
for all $x,y\in \Ln$. Let $\ad (z)(x)=\{z,x\}$ denote the adjoint operators of $\Ln$.
Taking $x=z$ in \eqref{post10} we obtain
\begin{align}\label{11}
\ad (z)^3+(2\la+1)\ad (z)^2+(\la^2+\la)\ad (z)& =0.
\end{align}
The minimal polynomial of $\ad (z)$ is a divisor of 
\[
t^3+(2\la+1)t^2+(\la^2+\la)t=t(t+\la+1)(t+\la).
\]
Hence the only possible eigenvalues of $\ad (z)$ are given by $0,-\la-1$ or $-\la$. \\[0.2cm]
{\it Case 1:} Assume that $\la=0$. Since $\Ln$ is semisimple it follows $\tr (\ad(z))=0$.
Hence all eigenvalues of $\ad(z)$ are zero and we have $\ad (z)^m=0$ for some $m\ge 1$.
Since \eqref{11} simplifies to $\ad(z)^3=-\ad (z)^2$ we obtain $\ad (z)^2=\cdots = \pm \ad(z)^m=0$.
This means that $z$ is a so called {\it sandwich element}. By Jacobi identity we have
$\ad(z)\ad(x)\ad(z)=0$ for all $x\in \Ln$. This implies $(\ad(z)\ad(x))^2=0$, so that
the Killing form of $\Ln$ satisfies $\tr (\ad(z)\ad(x))=0$ for all $x\in \Ln$. Since it
is non-degenerate it follows $z=0$. This implies $\phi=0$ and $x\cdot y=0$. \\[0.2cm]
{\it Case 2:} Assume that $\la=-1$. Because of $\tr (\ad(z))=0$ all eigenvalues of $\ad(z)$
are zero. The identity \eqref{11} reduces to $\ad(z)^3=\ad(z)^2$ and we obtain $z=0$ as before.
This implies $\phi=-\id$ and $x\cdot y=-\{x,y\}$. \\[0.2cm]
{\it Case 3:} Assume that $\la\neq 0,-1$. Now we may assume that $\ad(z)$ has three different eigenvalues
$0$,$-\la$, $-1-\la$, and its minimal polynomial is given by $t(t+\la+1)(t+\la)$ with
distinct factors. Hence $\ad(z)$ is diagonalizable and we are in the situation of Lemma $\ref{2.11}$. 
Therefore we have $0=\al+\be=-\la-1-\la$. This means $\la=-\frac{1}{2}$,  $\al=-\be=\frac{1}{2}$ and
$\Ln=V_0\oplus V_{\frac{1}{2}}\oplus V_{-\frac{1}{2}}$. For given $x,y\in V_0$ we have
$\{z,x\}=\{z,y\}=0$, and \eqref{post10} implies $(\la^2+\la)\{x,y\}=0$. By assumption we
can conclude that $\{x,y\}=0$. This shows that $\{V_0,V_0\}=0$. We also know that
$\{V_{\frac{1}{2}},V_{\frac{1}{2}}\}=\{V_{-\frac{1}{2}},V_{-\frac{1}{2}}\}=0$ 
from the proof of Lemma $\ref{2.11}$. Hence $\Ln$ is the direct vector space sum of three
abelian subalgebras $V_0,V_{\al},V_{\be}$. By Lemma $\ref{2.13}$ we have
$\Ln \simeq \Ls\Ll_2(\C) \oplus \cdots \oplus \Ls\Ll_2(\C)$, which we have excluded.
\end{proof}

One may ask what happens in the remaing case where $\Ln$ is a direct sum of 
copies of $\Ls\Ll_2(\C)$. It is enough to consider the case $\Ln=\Ls\Ll_2(\C)$. Let
$(e,f,h)$ be the standard basis with $\{e,f\}=h$, $\{h,e\}=2e$ and $\{h,f\}=-2f$.
Denote by $\Lr_{3,1}(\C)$ the solvable non-nilpotent Lie algebra with basis $e_1,e_2,e_3$ 
and brackets $[e_1,e_2]=e_2$, $[e_1,e_3]=e_3$.

\begin{prop}
Let $\Ln=\Ls\Ll_2(\C)$ and suppose that 
$x\cdot y =\{\{z,x\},y\}+\la \{x,y\}$ is a post-Lie algebra structure on $(\Lg,\Ln)$ for 
some $z\in \Ln$, $\la\in \C$. Then one of the following cases occurs.
\begin{itemize}
\item[(a)]  $z=0$, $\la=0$, $x\cdot y=0$ and $[x,y]=\{x,y\}$.
\item[(b)]  $z=0$, $\la=-1$, $x\cdot y=-\{x,y\}$ and $[x,y]=-\{x,y\}$.  
\item[(c)]  $z\neq 0$, $\la=-\frac{1}{2}$, $x\cdot y=\{\{z,x\},y\}-\frac{1}{2}\{x,y\}$
and $[x,y]=\{\{z,x\},y\}$. The Lie algebra $\Lg$ is isomorphic to $\Lr_{3,1}(\C)$.
\end{itemize}
\end{prop}

\begin{proof}
The arguments in the proof of Theorem $\ref{4.4}$ are also valid for $\Ln=\Ls\Ll_2(\C)$.
In the third case however we can classify all post-Lie algebra structures.
We first conclude $\la=-\frac{1}{2}$ and $x\cdot y=\{\{z,x\},y\}-\frac{1}{2}\{x,y\}$.
Then, writing $z=\al e+\be f+\ga h$, an easy calculation shows that this product  
is a post-Lie algebra structure on $(\Lg,\Ln)$ if and only if
\[
\frac{\al\be+\ga^2}{16}=1.
\]
In particular we have $z\neq 0$. The Lie brackets of $\Lg$ are given by
\begin{align*}
[e,f] & =-2\al e+2\be f, \\
[e,h] & =-4\ga e+2\be h,\\
[f,h] & =-4\ga f+2\al h.
\end{align*}
This is a solvable non-nilpotent Lie algebra which is isomorphic to $\Lr_{3,1}(\C)$
for all $\al,\be,\ga$ with $\al\be+\ga^2=1/16$.
\end{proof}

\end{document}